\pdfoutput=1 
\documentclass[11pt,oneside]{amsart}
\usepackage[leqno]{amsmath}
\usepackage{amssymb,amsthm} 
\usepackage[parfill]{parskip} 
\usepackage{libertine} 
\usepackage[T1]{fontenc}
\usepackage{textcomp}
\usepackage[varqu,varl]{inconsolata} 
\usepackage{amsmath,amsthm}
\usepackage[libertine,bigdelims,vvarbb]{newtxmath}
\usepackage[supstfm=libertinesups,%
  supscaled=1.2,%
  raised=-.13em]{superiors}
\usepackage{graphicx} 
\usepackage{calc}

\newcommand{\beq}{\begin{equation}}
\newcommand{\eeq}{\end{equation}}
\numberwithin{equation}{section} 
\newtheorem{thm}[equation]{Theorem}

\theoremstyle{remark}

\theoremstyle{definition}
\newtheorem{defn}[equation]{Definition}

\title{Subset Sum Instances in ZFC Limbo}
\usepackage[paperwidth=7in, paperheight=10in, textwidth=4.75in, textheight=8in, includehead=true]{geometry}
\usepackage{verbatim} 
\usepackage{url}
\usepackage[pdfborderstyle={/S/U/W 1}]{hyperref}

\hypersetup{
pdftitle={ZFC},
pdfauthor={\textcopyright\ S. G. Williamson}
}
\author{S. Gill Williamson}
\thanks{Department of Computer Science and Engineering, 
University of California San Diego; \url{http://cseweb.ucsd.edu/~gill/}; 
gwilliamson@ucsd.edu;
{\bf Keywords:} combinatorics; lattice graphs;  ZFC independence; P equals NP;
order type equivalence;  regressive regularity; subset sum problem
}

\begin{document}
\begin{abstract}
Our main result, Theorem~\ref{thm:polyprop}, proves the existence an infinity of subset sum problems solvable in polynomial time. The only proof we have of this result uses the ZFC independent Jump Free Theorem of Friedman~\cite{hf:alc}, thus putting Theorem~\ref{thm:polyprop} in what we informally call ZFC limbo~\cite{gw:lim}. The mathematics we use is elementary and at the level of a good undergraduate course in combinatorics or design and analysis of algorithms. The statement of the Jump Free Theorem is also easily understood at this level (but not the proof).  
\end{abstract}

 \maketitle
 
\section{Introduction}
Let $N=\{0,1,2,\ldots\}$ and $k\geq 2$.  In this note we derive some connections between the space and time hierarchies of computer science and the large cardinal hierarchies of set theory. Our main result, Theorem~\ref{thm:polyprop}, concerns families $S$ of functions $f$ with domains finite  subsets of $N^k$ and ranges finite subsets of $N$ ("finite" functions for short).  We define these families  $S$ below in terms of elementary set theoretic properties: finite, reflexive, jump free, full, etc.  

Let $k,p \geq 2$. For any pair $(f, E)$ and  $E^k\subseteq {\it domain} (f)$, $|E|=p$,  the function $f\in S$ naturally partitions $E^k\subset N^k$ into three blocks $E_0^k$, $E_1^k$, $E_2^k$ (Definition~\ref{def:fncpartitions}).
We define (Definition~\ref{def:fncpartsets}) sets of integers 
$\Delta[f]  E_0^k$, $\Delta[f]  E_1^k$,  $\Delta[f]  E_2^k$ associated with these three blocks. This construction is parameterized by functions $\rho$ and $\gamma$ (Definition~\ref{def:two ff}) and associated with integers $t\geq1$ (Definition~\ref{def:tlogbn})
so the sets $H_p=\Delta[f]  E_0^k\cup \Delta[f] E_1^k \cup \Delta[f]  E_2^k$ are "generic" in some sense.   We consider sequences of sets of this type
$H_2, H_3, \ldots, H_p, \ldots $. Note that $|H_p|=p^k$ (as a multiset). 
We prove 
Theorem \ref{thm:polyprop} that for fixed $t\geq 1,\,k\geq 2$ associated "jump free family" $S$,  there is at least one sequence of the type $H_2, H_3, \ldots, H_p, \ldots $ that, considered as instances to subset sum target $0$, can be solved in polynomial time $O((p^k)^t)$ where $p^k$ measures the size of the instances. Our only known proof uses Friedman's Jump Free Theorem~\ref{thm:jumpfree} which is independent of ZFC. We conclude with some interesting consequences.
\section{Basic definitions and theorems}

We denote by $N$ the set of all nonnegative integers. 
For $z=(n_1, \ldots, n_k)\in N^k$, $\max\{n_i\mid i=1,\ldots, k\}$ is
denoted by $\max(z)$.  Define $\min(z)$ similarly. 

In this section and the next we discuss Friedman's Jump Free Theorem.  In Section~3  we prove our main results.  
\begin{defn}[\bf Cubes and Cartesian powers in $N^k$]
The set  $E_1\times\cdots\times E_k$, where $E_i\subset N$, $|E_i|=p$, $i=1,\ldots, k,$   is called a $k$-cube of length $p$.  If $E_i = E, i=1,\ldots, k,$ then this cube is  $E^k$, the $k$th Cartesian power of $E$.
\label{def:cubespowers}
\end{defn}

\begin{defn}[\bfseries Equivalent ordered $k$-tuples]
\label{def:ordtypeqv}
Two k-tuples in $N^k$, $x=(n_1,\ldots,n_k)$ and $y=(m_1,\ldots,m_k)$, are  
{\em order equivalent tuples $(x\, ot \,y)$} if the sets
$\{(i,j)\mid n_i < n_j\} =  \{(i,j)\mid m_i < m_j\}$ and  $\{(i,j)\mid n_i = n_j\} =  \{(i,j)\mid m_i =m_j\}.$  
\end{defn}

Note that $ot$ is  an equivalence relation on $N^k$.
The number of equivalence classes is $\sum_{j=1}^k \sigma(k, j)< k^k,$ $k\geq 2$,
where $\sigma(k,j)$ is the number of surjections from a $k$ set to  a $j$ set.
We use ``$x\,ot\,y$'' and ``$x,\,y$ of order type $ot$'' to mean $x$ and $y$ belong to the same order type equivalence class.

We present some basic definitions due to Friedman~\cite{hf:alc}, \cite{hf:nlc}.

\begin{defn}[\bf Regressive value]
\label{def:regval}
Let $Y\subseteq N$, $X\subseteq N^k$ and $f:X\rightarrow Y$.  An integer $n$
is a  {\em regressive value} of $f$ on $X$
if there exist $x$ such that $f(x)=n<\min(x)$ .
\end{defn}

\begin{defn}[\bf Field of a function and reflexive functions]
\label{def:fldreflfncs}
For $A\subseteq N^k$ define ${\rm field}(A)$ to be the set of all coordinates of elements of $A$.  A function $f$ is {\em reflexive} in $N^k$ if 
${\rm domain}(f) \subseteq N^k$ and  ${\rm range}(f) \subseteq {\rm field}({\rm domain}(f))$.
\end{defn}

\begin{defn}[{\bf The set of functions} $T(k)$ ]
\label{def:tk}
$T(k)$ denotes all reflexive functions with finite domain.  We sometimes denote a function with domain $D$ by $f_D$
\end{defn}

\begin{defn} [\bf Full and jump free families]
\label{def:fulljump} 
Let $Q\subseteq T(k).$
\begin{enumerate}

\item {\bf full family:}  We say that $Q$ is a {\em full} family of functions on $N^k$ if for every finite subset 
$D\subset N^k$ there is at least one function $f$ in $Q$ whose domain is $D$.

\item{\bf jump free family:} For $D\subset N^k$ and $x\in D$ define $D_x = \{z\mid z\in D,\, \max(z) < \max(x)\}$. 
Suppose that for all $f_A$ and $f_B$  in $Q$, where $f_A$ has domain $A$ and $f_B$ has domain $B$,  the conditions
 $x\in A\cap B$, $A_x \subseteq B_x$, and $f_A(y) = f_B(y)$ for all $y\in A_x$ imply that 
$f_A(x) \geq f_B(x)$.  Then $Q$ will be called a {\em jump free} family of functions on $N^k$.  
\end{enumerate}
\label{def:fullrefjf}
\end{defn} 

\begin{figure}[h]
\begin{center}
\includegraphics[scale=.85]{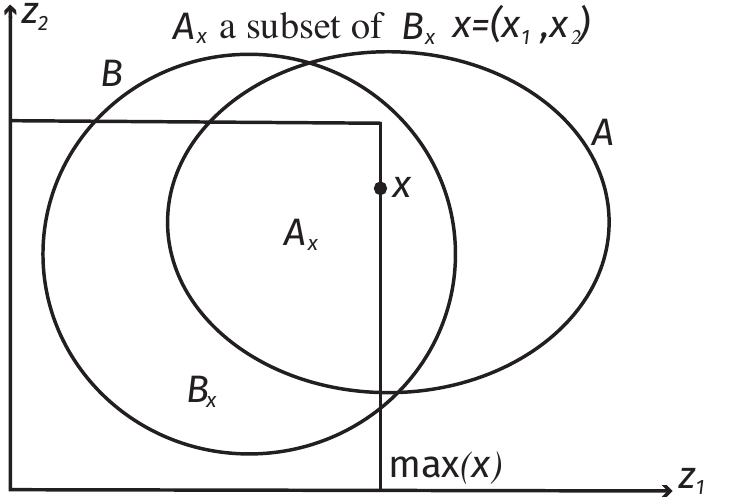}
\caption{Basic jump free condition~\ref{def:fullrefjf}.}
\label{fig:jfvenn}
\end{center}
\end{figure}

The jump free condition is illustrated in Figure \ref{fig:jfvenn}. Intuitively, referring to Figure~\ref{fig:jfvenn}, suppose that the region $A_x$ is to be searched for the smallest of some quantity and the result recorded at $x$.  Next, the search region is expanded to a superset $B_x$ with the search results for $A_x$ still valid 
(i.e., $f_A(y) = f_B(y)$ for all $y\in A_x$).  Then, clearly $f_A(x) \geq f_B(x)$.   

\begin{defn}[\bf Function regressively regular over $E$]
\label{def:regreg}
Let $k\geq 2$, $D\subset N^k$, $D$ finite, $f: D\rightarrow N$. 
We say that $f$ is {\em regressively regular} over 
$E$ (alternatively $f$ is {\em regressively regular over $E^k$)}  
$E^k\subset D$, if for each order type equivalence class $ot\,$ of $k$-tuples of $E^k$ either (1) or (2) occurs:
\begin{enumerate}
\item{\bf constant less than min $\mathbf{E}$:}  For all $x, y\,\in E^k$ of order type $ot$, $f(x)=f(y)< \min(E).$ 
\item{\bf greater or equal min:}   For all $x\in E^k$ of order type $ot$ $f(x)\geq \min(x).$
\end{enumerate}
\end{defn}

We use ZFC for the axioms of set theory: Zermelo-Frankel  plus the axiom of choice. 
The Jump Free Theorem, Theorem~\ref{thm:jumpfree} stated below, can be proved in 
 ZFC + ($\forall n$)($\exists$ $n$-subtle cardinal)
 but not in ZFC+ ($\exists$ $n$-subtle cardinal) for any fixed $n$ (assuming this theory is consistent).
A proof of the Jump Free Theorem is in Section 2 of \cite{hf:alc},
``Applications of Large Cardinals to Graph Theory.''  This proof references certain results from~\cite{hf:nlc}.  Most students can just take the proof of the Jump Free Theorem as given for this discussion. The proof is complicated (but well checked by logicians). Wikpedia is a good source for learning some basics about ZFC and large cardinals.

\begin{thm}[\bf Jump Free Theorem (\cite{hf:alc}, \cite{hf:nlc})] 
\label{thm:jumpfree}
Let $p, k\geq 2$ and $S\subseteq T(k)$ be  a full and jump free family of functions on $N^k$.
Then some $f_D \in S$ has at most $k^k$ regressive values on some 
$E^k \subseteq D$, $|E| = p$.  
In fact, some $f\in S$ is regressively regular over some $E$ of cardinality $p$.
\end{thm}
We note that the statement ``$f\in S$ is regressively regular over some $E$ of cardinality $p$'' implies that  ''some $f\in S$ has at most $k^k$ regressive values on some $E^k \subseteq {\rm domain}(f)$, $|E| = p$''.

\section{Subset sum instances and the Jump Free Theorem}
Let $S\subseteq T(k)$ be  a full and jump free family of functions on $N^k$.
The Jump Free Theorem asserts that for each $p\geq 2$ there exists an $f\in S$ such that $f$ is regressively regular over some $E^k\subseteq \mathrm{domain}(f)$ and $|E|=p$. 
Definition \ref{def:two ff} introduces two families of functions that will be used to construct instances to the subset sum problem. $Z$ denotes the integers.

\begin{defn}[\bf Families of functions] 
\label{def:two ff}
For each fixed $k\geq 2$ we define a two families of functions, $\rho$ and $\gamma$ with domains $E^k$, $E$ a finite subset of $N$.
\[
\mathbb{R}_k=\{\rho \mid \rho: E^k \rightarrow N,\,E^k\subset N^k,\, \rho(x)\geq \min(x),\, x\in E^k\}.
\]
\[
\mathbb{I}_k=\{\gamma\mid \gamma: E^k \rightarrow Z, E^k\subset N^k\}.
\]
 \end{defn}

In the next definition we consider functions $f$ and sets of the form $E^k$ contained in the domain of $f$.  We use $f$ to partition $E^k$  into three blocks which we call $E^k_0$, $E^k_1$, and $E^k_2$.  These partitions will be used to make a connection between sets of instances defined in Definition~\ref{def:fncpartsets} and regressive regularity, Definition~\ref{def:regreg}.

\begin{defn}[\bf $f\in S$ partitions $E^k$]
\label{def:fncpartitions}
Let $k, p\geq 2$, $f\in S$ , $E^k\subseteq {\rm domain}(f)$, 
$|E| = p.$  We define the blocks $E^k_0$,  $E^k_1$,  $E^k_2$ of a partition of $E^k$ as follows:
$$E^k_0 =\{ x\in E^k: f(x) < \min(E)\}$$
$$E^k_1=\{x\in E^k: \min(E)\leq f(x) < \min(x)\}$$
$$E^k_2=\{x\in E^k: \min(x)\leq  f(x) \}.$$
\end{defn}

We map the blocks of the partition~\ref{def:fncpartitions} into sets of integers as follows:
 
\begin{defn}[\bf Sets of instances from $E^k$]
\label{def:fncpartsets}
Let $k, p\geq 2$, $f\in S$ , $E^k\subseteq {\rm domain}(f)$,
Let the $\gamma$ and the $\rho$ be as in Definition~\ref{def:two ff}. Using Definition~\ref{def:fncpartitions} we 
 define $\Delta[f] E_0^k = \{f(x) - \min(E) : x\in E_0^k\}$,
$\Delta[f]  E_1^k = \{\gamma(x) : x\in E_1^k\}$ and $\Delta[f]  E_2^k = \{\rho(x)-\min(x) : x\in E^k_2 \}$. 
Let $H_p$ denote sets of the form 
$$H_p=\Delta[f]  E_0^k\cup \Delta[f] E_1^k \cup \Delta[f]  E_2^k $$
where $|E| = p$.
\end{defn}
We define a special class of functions $\rho$ of Definition~\ref{def:two ff}.

\begin{defn}[\bf $t$-log bounded $\rho$ over $E^k$, $|E|=p$]
\label{def:tlogbn} 
For $t\geq 1$ the function $\rho$ is $t$-log {\em bounded} over $E^k$ where $E=\{e_0, \ldots, e_{p-1}\}\,$  if  
$$
|\{\rho(x) - \min(x) : 0< \rho(x) - \min(x) < e_0 k^k,\, x\in E^k\}| \leq t\log_2 (p^k).
$$
\end{defn} 
Recalling that $\rho(x)\geq \min(x)$ and $\rho(x)$ can be arbitrarily large, we can choose 
$$
|\{\rho(x) - \min(x) :  \rho(x) - \min(x) \geq e_0 k^k,\, x\in E^k\}| 
$$
large enough to make $\rho$ $t$-log bounded.

\begin{thm}[\bf Instances to subset sum target zero]
\label{thm:polyprop}
Let $t\geq 1, k\geq 2$ be fixed. Let $p \geq 2$ and let $S\subseteq T(k)$ be  a full, reflexive and jump free family.
Referring to Definition~\ref{def:fncpartsets},  consider sets of the form $H_p=\Delta[f] E_0^k\cup \Delta[f] E_1^k \cup \Delta [f] E_2^k$ where $|E|=p$ and $f\in S$. We regard sequences of the form $H_2, H_3, \ldots , H_p, \ldots $, as  sequences of instances (size measured by $\,p^k$) to the subset sum problem target zero. For each fixed $t\geq 1$, $k\geq 2$, $S\subseteq T(k)$ there is a  sequence $H_2, H_3, \ldots , H_p, \ldots $ that can be solved in polynomial time where $p^k$ is the size of the instances.
\begin{proof}
For each $p \geq 2$ use the Jump Free Theorem~\ref{thm:jumpfree} to choose $f\in S$  regressively regular over some $E^k \subseteq {\mathrm domain(f)}$, $|E| = p$.

By regressive regularity, $\Delta[f] E_1^k$ is empty. There are at most $k^k$ terms in $\Delta[f] E^k_0$ and all are negative with absolute value less than 
$e_0 = \min(E)$. The terms in $\Delta[f] E^k_2$ are all greater than or equal to zero.  Zero is the target, so we assume all terms in
$\Delta[f] E^k_2$ are positive.  Since $E^k_2 \subseteq E^k$ we have
\[
|\{\rho(x) - \min(x) : 0< \rho(x) - \min(x) < e_0 k^k,\, x\in E_2^k\}| \leq t\log_2 (p^k).
\] 
Thus, to solve the subset sum problem we  compare all possible sums of $2^{k^k}$ subsets of negative terms with all $2^{t\log_2(p^k)}$ sums of subsets of positive terms.  
Positive terms greater than or equal $e_0 k^k$ are larger than any sums of negative terms.  Thus, we have $2^{k^k}(p^k)^t$ total comparisons of sums. Thus, for each fixed $t$ and $k$ the sequence of instances to  subset sum target zero, $H_2, H_3, \ldots , H_p, \ldots $ can be solved in polynomial time where $p^k$ is the size of the instances.
\end{proof} 
\end{thm}

\section{Conclusions} 
Theorem~\ref{thm:polyprop} uses the ZFC independent Jump Free Theorem (for each fixed  $t\geq 1$).
We know of no other proof thus putting  Theorem~\ref{thm:polyprop}  in what we call ''ZFC Limbo.'' 
If a ZFC proof could be found that the subset sum problem is solvable in polynomial time $O(n^\tau)$ where $n$ is the length of the instance ($p^k$ for fixed $k$ here), then that result would prove  Theorem~\ref{thm:polyprop} for $t=\tau$ and thus show it is provable within ZFC (i.e. remove it from ZFC limbo). 
We conjecture, however, that Theorem~\ref{thm:polyprop}  cannot be proved in  ZFC alone (for each fixed $t$). 
The basis for this conjecture is that the subset sum problem arises from
Theorem~\ref{thm:polyprop} and thus from the Jump Free Theorem in a very natural, generic way.
In particular, the sets of instances~\ref{def:fncpartsets} are constructed to be an alternative description of regressive regularity.
Of course, if our conjecture is true, ``subset sum is solvable in polynomial time'' cannot be proved in ZFC (perhaps because it is false).  

{\bf Acknowledgments:}  The author thanks  Professor Sam Buss (University of California San Diego, Department of Mathematics), 
Professor Emeritus Rod Canfield (University of Georgia, Department of Computer Science), and Professor Emeritus Victor W. Marek (University of Kentucky, Department of Computer Science) for numerous helpful suggestions.

%

\bibliographystyle{alpha}
\bibliography{jeffjoc}

\begin{thebibliography}{Wil17b}

\bibitem[Fri97]{hf:alc}
Harvey Friedman.
\newblock Applications of large cardinals to graph theory.
\newblock Technical report, Department of Mathematics, Ohio State University,
  1997.

\bibitem[Fri98]{hf:nlc}
Harvey Friedman.
\newblock Finite functions and the necessary use of large cardinals.
\newblock {\em Ann. of Math.}, 148:803--893, 1998.

\bibitem[RW99]{jg:pos}
Jeffrey~B. Remmel and S.~Gill Williamson.
\newblock Large-scale regularities of lattice embeddings of posets.
\newblock {\em Order}, 16:245--260, 1999.

\bibitem[Wil17a]{gw:lem}
S.~Gill Williamson.
\newblock Lattice exit models.
\newblock {\em arXiv:1907.11707v2 [math.CO]}, 2017.

\bibitem[Wil17b]{gw:sub}
S.~Gill Williamson.
\newblock On the difficulty of proving P=NP in Z{F}{C}.
\newblock {\em arXiv:1708.08186 [math.CO]}, 2019.

\bibitem[Wil17c]{gw:lim}
S.~Gill Williamson.
\newblock Combinatorics in ZFC limbo.
\newblock {\em JOC Vol. 10, Num. 3, 579-593} 2019.

\end{thebibliography}

\end{document}